\documentclass{amsart}
\newtheorem{thm}{Theorem}[section]

\newtheorem{prop}[thm]{Proposition}
\theoremstyle{definition}

\theoremstyle{remark}

\numberwithin{equation}{section}
\def\ed{\mathop{\mathrm{ext-dim}}}

\begin{document}

\title[Topological semigroups and universal
spaces]{Topological semigroups and universal spaces related to
extension dimension}
\author{A.~Chigogidze}
\address{Department of Mathematics and Statistics,
University of Saskatche\-wan, McLean Hall, 106 Wiggins Road,
Saskatoon, SK, S7N 5E6, Canada}
\email{chigogid@math.usask.ca}
\thanks{The first author was partially supported by NSERC research grant.}

\author{A.~Karasev}
\address{Department of Mathematics and Statistics,
University of Saskatche\-wan, McLean Hall, 106 Wiggins Road,
Saskatoon, SK, S7N 5E6, Canada}
\email{karasev@math.usask.ca}

\author{M.~Zarichnyi}
\address{Department of Mechanics and Mathematics, Lviv State University,
Universitetska 1, 290602 Lviv, Ukraine}
\email{topos@franko.lviv.ua}
\subjclass{54H15, 54F45, 55M10}
\keywords{Extension dimension, universal space, semigroup}
\date{}
\dedicatory{}
\commby{}

%%% ----------------------------------------------------------------------

\begin{abstract} It is proved that there is no structure of left
(right) cancelative semigroup on $[L]$-dimensional universal space
for the class of separable compact spaces of extensional
dimension $\le[L]$. Besides, we note that the homeomorphism group of
$[L]$-dimensional space whose nonempty open sets are universal for 
the class of separable compact spaces of extensional
dimension $\le[L]$ is totally disconnected.
\end{abstract}

%%% ----------------------------------------------------------------------
\maketitle
%%% ----------------------------------------------------------------------
\section{Preliminaries}
Let $L$ be a CW-complex and $X$ a Tychonov space. The {\it
Kuratowski notation\/} $X\tau L$ means that, for any continuous map
$f\colon A\to L$ defined on a closed subset $A$ of $X$, there
exists an~extension $\bar f\colon X\to L$ onto X. This notation
allows us to define the preorder relation $\preceq$ onto the class
of CW-complexes: $L\preceq L'$ iff, for every Tychonov space $X$,
$X\tau L$ implies $X\tau L'$ \cite{Ch}.

The preorder relation $\preceq$ naturally generates the equivalence
relation $\sim$: $L\sim L'$ iff $L\preceq L'$ and $L'\preceq L$. We
denote by $[L]$ the equivalence class of $L$.

The following notion is introduced by A.~Dranishnikov (see, \cite{D} and \cite{DD}). The
{\it extension dimension\/} of a~Tychonov space $X$ is less than or
equal to $[L]$ (briefly, $\ed(X)\le[L]$) if $X\tau L$.

We say that a Tychonov space $Y$ is said to be a~universal space for 
the class of compact
metric spaces $X$ with $\ed(X)\le[L]$ if $Y$ contains a topological copy 
of every compact
metric space $X$ with $\ed(X)\le[L]$. See \cite{Ch} and \cite{ChV} for 
existence of universal spaces.

In what follows we will need the following statement which appears in \cite{ChZ} as Lemma 3.2.
\begin{prop}\label{p:1}
Let $i_0=\min\{i:\pi_i(L)\neq0\}$. Then $\ed(S^{i_0})\le[L]$.
\end{prop}

\section{Main theorem}
Recall that a semigroup $S$ (whose operation is denoted as multiplication) 
is called a {\it left cancelation
semigroup\/} if $xy=xz$ implies $y=z$ for every $x,y,z\in S$.

\begin{thm} Let $L$ be a connected CW-complex and Let $Y$ be a~universal space for the class of compact
metric spaces $X$ with $\ed(X)\le[L]$. If $\ed(Y)=[L]$, then there
is no structure of left (right) cancelation semigroup on $Y$
compatible with its topology.
\end{thm}
\begin{proof}
Suppose the contrary and let $Y$ be a left cancelation semigroup.
Let  $\alpha(\coprod_{j=1}^\infty S^{i_0}_j)$ be the Alexandrov
compactification of the countable topological sum of copies $S^{i_0}_j$ of 
the sphere
$S^{i_0}$, where $i_0=\min\{i:\pi_i(L)\neq0\}$. By the countable sum 
theorem for extension dimension and  Proposition \ref{p:1},  $\ed(
\alpha(\coprod_{j=1}^\infty S^{i_0}_j))\le[L]$ and, since $Y$ is universal, 
$Y$ contains a~copy of
 $\alpha(\coprod_{j=1}^\infty S^{i_0}_j)$. We will assume that
$\alpha(\coprod_{j=1}^\infty S^{i_0}_j)\subset Y$. Besides, since
$\ed(Y)\ge[S^1]$, we see that $Y$ contains an arc $J$. Let $a,b$ be
endpoints of $J$. There exists $j_0$ such that $aS^{i_0}_{j_0}\cap
bS^{i_0}_{j_0}=\emptyset$. By Proposition \ref{p:1}, there exists
a~map $f\colon aS^{i_0}_{j_0}\cup bS^{i_0}_{j_0}\to L$ such that
$f|aS^{i_0}_{j_0}$ is a~constant map and $f|bS^{i_0}_{j_0}$ is not
null-homotopic. Extend map $f$ to a map $\bar f\colon Y\to L$. Let
$g\colon [0,1]\to J$ be a homeomorphism, then the map $F\colon
S^{i_0}_{j_0}\times[0,1]\to L$, $\bar f(x,t)=g(t)x$, is a homotopy
that contradicts to the fact that $f|bS^{i_0}_{j_0}$ is not
null-homotopic.
\end{proof}

The homeomorphism group $\mathrm{Homeo}(X)$ of a space $X$ is
endowed with the compact-open topology.

\begin{thm}\label{t:2}
Suppose $\ed(X)=[L]$ and every nonempty open subset of $X$ is
universal for the class of separable metric spaces $X$ with
$\ed(X)\le[L]$. Then the homeomorphism group $\mathrm{Homeo}(X)$ is
totally disconnected.
\end{thm}
\begin{proof}
Suppose the contrary. Let $h\in\mathrm{Homeo}(X)$,
$h\neq\mathrm{id}_X$. There exists $x\in X$ such that $h(x)\neq x$
and, therefore, there exists a neighborhood $U$ of $x$ such that
$h(U)\cap U=\emptyset$. Since $U$ is universal for the class of
separable metric spaces $X$ with $\ed(X)\le[L]$, there exists an
embedding of $S^{i_0}$ into $U$, where $i_0$ is as in Proposition
\ref{p:1}. We may suppose that $S^{i_0}\subset U$. There exists a
map $f\colon S^{i_0}\cup h(S^{i_0})\to L$ such that the restriction
$f|S^{i_0}$ is not null-homotopic while the restriction
$f|h(S^{i_0})$ is null-homotopic. Since $\ed(X)\le[L]$, there exists an
extension $\bar f\colon
X\to L$ of the map  $f$. The set $$W=\{g\in \mathrm{Homeo}(X): \bar
f|g(S^{i_0})\text{ is not null-homotopic }\}$$ is an open and
closed subset of $\mathrm{Homeo}(X)$. We see that $W$ is a neighborhood 
of unity that does not contain $h$.
\end{proof}

\section{Open problems}

Note that the case $L=S^n$ corresponds to the case of covering dimension. 
In this case, the topology of homeomorphism groups of some universal spaces has been
investigated  by many authors (see the survey \cite{CKT}).

In particular, it is known (see \cite{B} and \cite{OT}) that the homeomorphism 
group of the $n$-dimensional Menger compactum $M^n$ (note that $M^n$ satisfies the conditions of 
Theorem \ref{t:2} with $L=S^n$) is one-dimensional.

Let $[L]\ge[S^1]$ and $X$ be as in Theorem \ref{t:2}. Is
$\dim(\mathrm{Homeo}(X))\ge1$?

Another version: Is there $X$ that satisfies the conditions of 
Theorem \ref{t:2} and such that $\dim(\mathrm{Homeo}(X))\ge1$?

% ------------------------------------------------------------------------

\end{document}